\documentclass[a4,12pt]{article}
\usepackage{amsfonts,amssymb,mathrsfs,graphicx,amsbsy,color}

\def\d{{\rm d}}
\def\mi{{\rm i}}
\def\eps{\varepsilon}
\def\g{\gamma}
\def\G{\Gamma}
\def\l{\lambda}
\def\s{\sigma}
\def\z{\zeta}
\def\ZA{\mathop{\mathscr Z}\nolimits}
\def\ZB{\mathop{\mathcal Z}\nolimits}
\def\ZC{\mathop{\mathfrak Z}\nolimits}
\def\ZS{\mathop{\mathbf Z}\nolimits}
\def\Re{\mathop{\rm Re\,}\nolimits}
\def\Im{\mathop{\rm Im\,}\nolimits}
\def\e{\mathop{\rm e}\nolimits}
\def\Res{\mathop{\rm Res}\nolimits}
\def\hf{{\textstyle{1 \over 2}}}
\def\qt{{\textstyle{1 \over 4}}}

\def\defi{\stackrel{\rm def}{=}}
\def\si{\!\!\! &}
\def\se{& \!\!\!}
\def\sep{\,|\,}
\def\ni{\noindent}

\definecolor{dgray}{rgb}{.4,.4,.4}

\newcommand{\beq}{\begin{equation}}
\newcommand{\eeq}{\end{equation}}
\newcommand{\bea}{\begin{eqnarray}}
\newcommand{\eea}{\end{eqnarray}}
\newcommand{\tcA}{\textsf{\textbf{A}}}
\newcommand{\tcF}{\textrm{\textbf{F}}}
\newcommand{\tcR}{\textsf{\textbf{R}}}
\newcommand{\tcT}{\textit{\textbf{T}}}

\title{
Zeta functions \emph{over zeros of} Zeta functions\\
and an exponential-asymptotic view\\
of the Riemann Hypothesis}
\author{{\bf Andr\'e Voros}\\
Institut de Physique Th\'eorique, CEA-Saclay (CNRS-URA 2306)\\
F-91191 Gif-sur-Yvette Cedex (France)\\
e-mail: {\tt andre.voros@cea.fr}}

\begin{document}
\date{\emph{dedicated to Professor Takashi AOKI for his 60\textsuperscript{\,th} birthday}}
\maketitle

\begin{abstract}
We review generalized zeta functions built over the Riemann zeros (in short: ``superzeta" functions).
They are symmetric functions of the zeros that display a wealth of explicit properties,
fully matching the much more elementary Hurwitz zeta function. 
As a concrete application, a superzeta function enters
an integral representation for the Keiper--Li coefficients,
whose large-order behavior thereby becomes computable by the method of steepest descents;
then the dominant saddle-point entirely depends on the Riemann Hypothesis being true or not,
and the outcome is a sharp exponential-asymptotic criterion for the Riemann Hypothesis
that only refers to the large-order Keiper--Li coefficients. 
As a new result, that criterion, then Li's criterion, are transposed 
to a novel sequence of Riemann-zeta expansion coefficients based at the point 1/2 
(vs 1 for Keiper--Li).
\end{abstract}

\bigskip

It is a great honor and pleasure to dedicate this work to Professor AOKI, 
who has had zeta functions among his numerous activities; 
we are specially grateful to him for setting up a major conference in Osaka in 2003 
\cite{AKNO} that much stimulated our work afterwards.
And the present work also touches another theme in which we have much longer 
fruitfully interacted with him (and collaborators, mainly at the RIMS),
namely complex exponential asymptotics (here embodied in the method of steepest descents).

\section{(Generalized) zeta functions}

Discrete numerical sets $\{ w_k \}$ that have ``natural", ``collective" definitions
(e.g., roots of equations, spectra, \ldots),
are often better accessible through their \emph{symmetric functions},
which easily display richer and more explicit properties than the initial data.
Specially fruitful symmetric functions turn out to be the \emph{zeta functions}, defined (formally) as
\beq
\label{ZZD}
Z(x) = \sum_k w_k^{\, -x} \ \ \mbox{(standard zeta)}, \qquad 
Z(x,w) = \sum_k (w_k+w)^{-x} \ \ \mbox{(generalized zeta)},
\eeq
typically for complex arguments. As a rule, the exponent (here, $x$) is the key variable; 
the shift parameter ($w$) adds useful flexibility but retains an auxiliary status.

Based on \cite{VB} and references therein, we will specifically review cases 
where $\{ w_k \}$ is the set of nontrivial zeros $\{ \rho \}$ of Riemann's zeta function 
(or some related set, or a generalization) -
we will name the resulting functions (\ref{ZZD}) ``superzeta" functions for brevity. 

Our starting set $\{ w_k \}$ will however be the natural integers, with the zeta functions
\beq
\label{RHD}
\z (x) = \sum\limits_{k=1}^\infty k^{-x} \mbox{ (Riemann)}, \quad 
\z (x,w) = \sum\limits_{k=0}^\infty (k+w)^{-x} \mbox{ (Hurwitz)} \qquad (\Re x >1) 
\eeq
(the standard notations, in which $\z (x) \equiv \z (x,w=1)$), because

a) $\z (x)$ gives rise to the first superzeta functions, those over the \emph{Riemann zeros};

b) $\z (x,w)$ provides the basic \emph{template}
for the full set of explicit properties ultimately displayed by those less elementary
superzeta functions.

\subsection{Predicted pattern of general results for zeta functions}
\label{SPZ}

For our later sets $\{ w_k \}$ (all countably \emph{infinite}), typical results can be categorized as:
\medskip

\sffamily
\ni \textbullet \ \tcA\ - \textbf{Analytic structure in the whole complex $x$-plane} (at fixed $w$): 
our zeta functions will be defined by (\ref{ZZD}) only over a half-plane in $x$, as in (\ref{RHD}), 
but will then admit a \emph{meromorphic} continuation to the whole $x$-plane,
with \emph{exactly and fully computable poles and principal parts}. 
(E.g., Hurwitz's $ \z (x,\cdot)= 1/(x-1) + $ [entire function].)
\medskip

\normalfont
\ni \textbullet \ \tcF\ - {\bf Functional relation} (often): 
this may link the analytic continuation of $Z(x,w)$ 
to another function of Mellin-transform type (as in (\ref{JR}) for the Hurwitz case);
when it exists, \tcF\ may also readily supply all the explicit features announced here.
\medskip

\ni $\bullet\ \frac{\mbox{\tcR}}{\overline{\mbox{\tcT}}}$ - {\bf Special-value formulae}
for $Z(x,w)$ at all integers ${x=\mp n}$:
$ \frac{\mbox{\textsf{\textbf{Rational}}}}{\overline{\mbox{\textit{\textbf{Transcendental}}}}} $
respectively (such will be the layout of the Tables), and having quite distinct origins:
\vskip -1mm

\ni \rule{\textwidth}{.4pt}
\vskip -1mm

\sffamily
\ni \tcR\ \ for $x=-n$ (in conjunction with \tcA): 
our zeta functions relate to \emph{Mellin transforms}, e.g.,
${Z (x,w) = \G (x)^{-1} \int_0^\infty f(z,w) \, z^{x-1} \,\d z}$ (like (\ref{HM}) for Hurwitz)
where $f(z,w)$ will have an explicit $z \to 0$ power expansion:
this feature makes the Mellin integral meromorphic for all~$x$ by continuation
and also yields its explicit poles and principal parts (\tcA),
plus (upon the division by $\G (x)$) explicit \textbf{rational} (\tcR) values for $Z(x,w)$ 
at all $x=-n$ (0 \emph{and the negative integers}).
By ``rational" we mean: polynomial in $w$ and in the expansion coefficients of $f$ at $z=0$.
\vskip -3mm

\ni \rule{\textwidth}{.4pt}
\vskip -4.5mm
\ni \rule{\textwidth}{.4pt}
\vskip -.5mm

\normalfont\itshape
\ni \tcT\ \ \textit{for $x=+n:$ zeta functions relate to \emph{zeta-regularized} products
$D(w) = \prod\limits_k (w+w_k)$ by}
\vskip -5mm

\bea
\label{DET}
\log D(w) \si\defi\se - \partial_x Z(x,w)\vert_{x=0} \qquad \mbox{for }
Z(x,w) = \textstyle \sum\limits_k (w_k+w)^{-x} \\[-6pt]
\Longrightarrow \qquad Z(n,w) \si=\se \frac{(-1)^{n-1}}{(n-1)!} (\log D)^{(n)}(w) 
\ \mbox{for } n=1,2,\ldots \quad \mbox{outside poles of } Z ,
\eea
\vskip -5mm

\beq
\label{FPN}
\mbox{vs } \left\{ \matrix{ 
{\rm FP}_{x=1} Z(x,w) \si=\se (\log D)'(w) \quad \mbox{\emph{(FP: finite-part} extraction)} \cr
{\rm FP}_{x=n} Z(x,w) \si=\se \mbox{\emph{[a bulkier formula]}\, for $n \ge 2$ \emph{\cite[\S \,2.6.1]{VB}}} } \right\}
\mbox{ on poles of }Z \,
\eeq
\textit{(the very last case is unused here); thus, formulae (\ref{DET})--(\ref{FPN})
express zeta values at all $x=+n$ (\emph{0 and the positive integers}),
in terms of a function $D(w)$ which is entire, has $\{-w_k\}$ as zeros, 
and will be \textbf{transcendental} (\tcT) but hopefully known. 
E.g., in the Hurwitz case: $\{ w_k \} =\{ 0,1,2,\ldots \} $ leads to}
\beq
\label{GDet}
D(w) = \sqrt{2\pi} / \G (w).
\eeq
\vskip -4.5mm

\ni \rule{\textwidth}{.4pt}

\normalfont
\subsection{Results for the Hurwitz zeta function $\z (x,w)$}

We recall classic results \cite{BF}, to be needed for the notations and for later reference.
\medskip

\begin{figure}[h]
\centering
\includegraphics[scale=1.2]{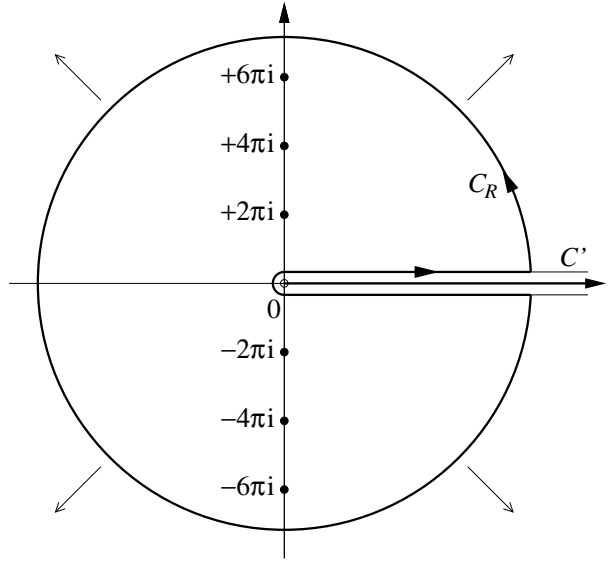}
\caption{\label{F1} \small
Integration paths in (\ref{HI})--(\ref{HR}) leading to the functional relation (\ref{JR}) 
for $\z (x,w)$.
}
\end{figure}

\ni \textbullet \ \tcF - Functional relation (will imply \tcA): we start from 
the (obvious) Mellin-transform representation of $\z (x,w)$, then convert it to a Hankel integral:
\bea
\label{HM}
\z (x,w) \si=\se \frac{1}{\G (x)} \int_0^\infty \frac{\e^{(1-w)z}}{\e^z-1} \, z^{x-1} \,\d z
\qquad \qquad \qquad \quad (\Re x > 1, \ \Re w>0) \\
\label{HI}
\si=\se \frac{\G (1-x)}{2\mi\pi} \int_{C'} \frac{\e^{(1-w)z}}{\e^z-1} \, (-z)^{x-1} \,\d z
\qquad \qquad (x\not\in {\mathbb N}^\ast , \ 1 \ge \Re w >0) \\
\label{HR}
\si=\se \frac{\G (1-x)}{2\mi\pi} \lim_{R \to +\infty} \int_{C_R} \frac{\e^{(1-w)z}}{\e^z-1} \, (-z)^{x-1} \,\d z
\quad (x\not\in {\mathbb N}^\ast , \ 1 \ge \Re w >0) 
\eea
(see Fig.~\ref{F1}); then this contour integral is readily evaluated by the residue calculus, giving
\bea
\label{JR}
\z (x,w) \si=\se \frac{\G (1-x)}{(2\pi )^{1-x} \,\mi}
\bigl[ \e^{\mi\pi x/2} F(\e^{2\mi\pi w},1-x) - \e^{-\mi\pi x/2} F(\e^{-2\mi\pi w},1-x) \bigr] 
\qquad \mbox{(\tcF)} \qquad \\
\label{LP}
\mbox{for }F(u,y) \si\defi\se \sum\limits_{n=1}^\infty \frac{u^n}{n^y}, \quad 
\mbox{and  } \Re x <0, \ 1 \ge \Re w >0 \qquad \qquad \qquad \qquad \quad \; \mbox{(\tcA)}
\eea
(Lerch or polylogarithm function; (\ref{JR}): Jonqui\`ere's formula).

\ni $\bullet\ \frac{\mbox{\tcR}}{\overline{\mbox{\tcT}}}$ - Special values: \tcR\ resulting from (\ref{HM}), 
followed by \tcT\ from (\ref{GDet}), are listed in Table~\ref{T1}. 

\begin{table}[h]
\centering
\begin{tabular}  {ccc}
\hline \\[-12pt]
$x$ & $ \z (x,w) = \sum\limits_{k=0}^\infty (k+w)^{-x} $ & \\[8pt]
\hline \\[-10pt]
$-n \le 0$ & $\mathsf{ -\frac{\textstyle 1}{\textstyle n + 1} B_{n+1}(w) }$ & \\[8pt]
$0$ & $\mathsf{ \hf - w }$ & \tcR \\[3pt]
\hline \\[-13pt]
\hline \\[-13pt]
$\textstyle 0 \atop \mbox{\it \small ($x$-derivative)}$ &
 ~~~~\fbox{$ \z '(0,w) = \log \bigl( \G (w)/\sqrt {2\pi } \bigr) $}~~~~ & \tcT \\[5pt]
$\textstyle +1~~ \atop \mbox{\it \small (finite part)}$ & 
\fbox{$ {\rm FP}_{x=1} \z (x,w) = - \bigl( \log \G \bigr)' (w)$} & \\[4pt]
$+n >1$ & $ \displaystyle \frac{(-1)^n}{(n-1)!} \, \bigl( \log \G \bigr) ^{(n)}(w)$ &
\\[6pt]
\hline
\end{tabular}
\caption{\label{T1} \small 
Special-value formulae for the Hurwitz zeta function,
$ \frac{\mbox{\textsf{Rational}}}{\overline{\mbox{\textit{Transcendental}}}} $ at all integer $x=\mp n$.
($\mathsf{ B_n(w): }$ \textsf{Bernoulli polynomials}.)
}
\end{table}

\medskip

Moreover, (\ref{JR}) reduces at $w=1$ to \emph{Riemann's Functional Equation}: 
\beq
\label{RFE}
\z (x) \equiv 2 (2\pi )^{x-1} \sin \frac{\pi x}{2} \, \G (1-x) \, \z (1-x) 
\quad \iff \quad \Xi (\hf+t) \equiv \Xi (\hf-t),
\eeq
in terms of a \emph{completed} zeta function
\beq
\label{CZ}
\Xi (x) = x(x-1) \pi^{-x/2} \G (x/2) \, \z (x) \quad (\mbox{our normalization: } \Xi (0) = \Xi (1) =1),
\eeq
and of the convenient variable $ t \equiv x-\hf$. Remarkable values of $t$ are then
\medskip

\centerline{$t=0$ (the center of symmetry for $\Xi (\hf + \cdot)$) \quad and 
$\quad t=+\hf$ (the pole of $\z (\hf + \cdot)$).}
\medskip

Due to (\ref{RFE}), $\z (x)$ has \emph{trivial zeros}: $x=-2k,\ k=1,2,\ldots$, 
plus those of $\Xi (x) :$ the \emph{Riemann zeros} $\rho $; 
the latter are countably many, 2-by-2 symmetrical about $t=0$, as
\beq
\rho = \hf \pm \mi \tau_k \qquad (\Re \tau_k >0) \qquad (k=1,2,\ldots),
\eeq
and lie within the \emph{critical strip} $\{ 0 < \Re \rho < 1 \}$; 
moreover, the \emph{Riemann Hypothesis} \cite{Ri} (still an open conjecture) 
puts them all on the \emph{critical line}:
\beq
\label{RRH}
\Re \rho = \hf \quad \Longleftrightarrow \quad \tau _k 
\mbox{ is real \qquad for all Riemann zeros. \qquad \qquad \bf [RH]}
\eeq

\section{Superzeta functions, in the Riemann zeros' case}

Over the set of Riemann zeros $\{ \rho \}$ (always counted with multiplicities, if any),
several generalized zeta functions are conceivable:
\bea
\label{ZAD}
\ZA (s \sep t) \si=\se \sum_\rho (\hf+t - \rho)^{-s}, \qquad \qquad
\Re s > 1 \quad \mbox{(1\textsuperscript{st} kind}) \\
\label{ZBD}
\ZB (\s \sep t) \si=\se \sum_{k=1}^\infty ({\tau_k}^2 + t^2)^{-\s }, \ \qquad \qquad
\Re \s > \hf \quad \mbox{(2\textsuperscript{nd} kind}) \\
\label{ZCD}
\ZC (s \sep \tau) \si=\se \sum_{k=1}^\infty (\tau_k + \tau)^{-s}, \qquad \qquad \quad \,
\Re s > 1 \quad \mbox{(3\textsuperscript{rd} kind}) 
\eea
where the parameter $t$ has the same meaning as in (\ref{RFE}); those functions are 2-by-2 inequivalent,
except at a single parameter location where a \emph{confluence} occurs:
\vskip -3mm

\beq
\label{Z0}
(2 \cos \hf \pi s)^{-1} \! \ZA (s \sep 0) \equiv \ZB(\hf s \sep 0) \equiv \ZC (s \sep 0)
\equiv \sum_{k=1}^\infty \tau _k^{\, -s} ;
\eeq
(\ref{Z0}) can be taken as a standard (one-variable) zeta function over the Riemann zeros.

Those zeta functions over the Riemann zeros (or ``superzeta" functions, for brevity) 
got considered quite sporadically until the turn of the century:
we found but a dozen studies, all partially focused and incomplete even as a whole (Table~\ref{T2}).
From 2000 onwards we have gradually filled gaps in the description \cite{Vz}\cite{VZ}\cite{VO},
ending up with a more global perspective, also set in new coherent notations, 
as a book \cite{VB}. The present text then basically condenses the main results of \cite{VB}, 
where the details skipped here can be found.

\begin{table}
\small \centering 
\begin{tabular}  {cccc}
\hline \\[-10pt]
year & $\ZA(s \sep t) = \sum\limits_\rho (\hf+t - \rho)^{-s}$ & 
$\!\!\!\!\!\!\!\!\! \ZB(\s \sep t) = \sum\limits_{k=1}^\infty (\tau_k^{\, 2} + t^2)^{-\s }$ & 
$\!\!\! \ZC(s \sep \tau) =\sum\limits_{k=1}^\infty (\tau_k + \tau)^{-s}$ \\[4pt]
\hline \\[-8pt]
\sl 1860 & Riemann (unpubl.): $\ \sum\limits_\rho \, \rho^{-1}$ \tcT && \\
1917 & Mellin \ $t=\pm \hf$ \tcA \tcR && $\!\!\!\!\!\!$ Mellin \ $\tau = \pm \hf \mi$ \,\tcA \\
$ \matrix{\rm 1949 \cr 1966 \cr 1970 \cr 1971} $ & 
\multicolumn{3}{c}{\qquad $ \hbox to 60pt{\leftarrowfill} \quad
\left. \matrix{\mbox{Guinand \tcA \tcF} \cr \mbox{Delsarte \tcA} \cr \mbox{Chakravarty\,\tcF}\cr \mbox{\qquad " \quad \ \tcA \tcR}} 
\!\! \right\} t=0 \mbox{ (confluent case)} \quad \hbox to 60pt{\rightarrowfill} $} \\
$ \matrix{\rm 1985 \cr 1988} $ &
$ \!\!\!\! \left. \matrix{\rm Matsuoka \cr \rm Lehmer} \!\! \right\} \sum\limits_\rho \, \rho^{-n} 
\left\{ \! \matrix{t=\hf \cr s = n} \right. $\tcT \\
1988 & & Kurokawa \qquad $t=\hf$ \tcA~~~~~~~ & \\
1989 & & $\!\!\!\!\!\!\!\!\!\!\!\!\!\!\!\!\!\!\!\!\!$ Matiyasevich 
$\! \left\{ \! \matrix{\, t= \hf \hfill \cr \s =n} \right. $\tcT \\
$ \matrix{\rm 1992 \cr 1994} $ & 
$ \!\!\!\! \left. \matrix{\rm Keiper \cr \mbox{Zhang--Williams}} \!\! \right\} \sum\limits_\rho \, \rho^{-n}$ \tcT \\
$ \matrix{\rm 1992 \cr 1994} $ &
$ \!\! \left. \matrix{\rm Deninger \cr \mbox{Schr\"oter--Soul\'e}} \!\! \right\}
\partial_s \! \ZA (0 \sep t) \left\{ \! \matrix{\mbox{\tcA} \cr \mbox{\tcT}} \right. $ \\[9pt]
\hline
\end{tabular}
\caption{\label{T2} \small
The pre-2000 literature on ``superzeta" functions of all 3 kinds (to our knowledge),
sorted by kind (columns). (\tcA--\tcT: classes of results, cf. Sect.~\ref{SPZ};
\ $n$ means any positive integer; the references are listed in full at the end.)
}
\end{table}

\subsection{The superzeta function of the 1\textsuperscript{st} kind $\ZA (s \sep t)$}

\subsubsection{Basic functional relation / continuation formula}

$\ZA (s \sep t)$ is only defined by (\ref{ZAD}) where its series converges,
i.e., in the half-plane $\{ \Re s >1 \}$.

Now, in terms of the partner zeta function over the \emph{trivial zeros} of $\z (x)$,
\beq
\label{sz1}
\ZS (s \sep t) \defi \sum_{k=1}^\infty (\hf+t+2k)^{-s} \equiv 
2^{-s}\z (s,{\textstyle{5 \over 4}} + \hf t) ,
\eeq
and of another \emph{Mellin transform},  
\beq
\label{ZRP}
{\mathscr J} (s \sep t) = \int\limits_0^\infty \frac{\z '}{\z } (\hf +t+y) \, y^{-s} \, \d y \qquad (\Re s <1) ,
\eeq
$\ZA (s \sep t)$ admits an analytic continuation formula which reads as
\beq
\label{ZC1}
\fbox{$\displaystyle \ZA (s \sep t) = -\ZS (s \sep t) + (t-\hf)^{-s}
+ \frac{\sin \pi s}{\pi} {\mathscr J} (s \sep t) $} \qquad (\Re s <1) ,
\eeq
valid in a cut $t$-plane: $\hf+t$ must avoid all negatively oriented half-lines 
drawn from the Riemann zeros and the pole (whose cut $t \in (-\infty,+\hf]$ imposes
special treatments for both remarkable values $t=0$ and $\hf$).

\begin{figure}[h]
\centering
\vskip -2mm
\includegraphics[scale=1.2]{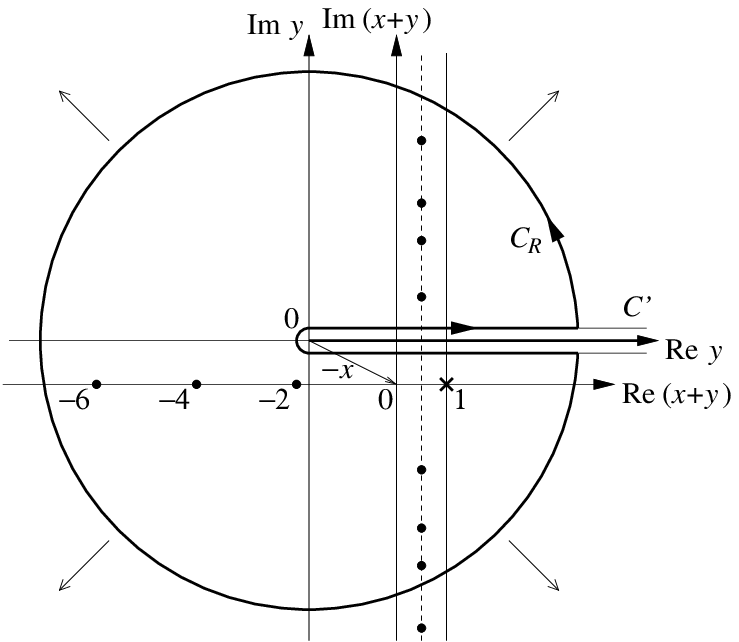}
\vskip -2mm
\caption{\label{F2} \small
Integration paths applied to the Mellin transform (\ref{ZRP}) similarly to Fig.~\ref{F1}, 
leading to the analytic continuation formula (\ref{ZC1}) for $\ZA (s \sep t)$. 
($x \equiv \hf+t$;
\protect\raisebox{1pt}{$\ \scriptstyle\bullet $} : zeros of $\z (x+y)$, 
including mock ``Riemann zeros" with imaginary parts contracted to fit into the frame.)
}
\end{figure}

The argument of proof (parallel to (\ref{HM})--(\ref{HR}), and sketched by Fig.~\ref{F2}) shows 
(\ref{ZC1}) to be the counterpart of the functional relation (\ref{JR}) from the Hurwitz case,
i.e., (\ref{ZC1}) gives the property {\tcF} for~$\ZA $.

Moreover, by the same logic as in Sect.~\ref{SPZ}.\tcR, 
the $y \to 0$ Taylor expansion of the integrand in (\ref{ZRP})
induces the continuation of $\mathscr J (s \sep t)$ (and of the \tcF-identity (\ref{ZC1})) 
to the whole $s$-plane, and it fully specifies the poles and residues of $\mathscr J$. 
As a result, the product $[\sin \pi s \mathscr J (s \sep t) ]$ in (\ref{ZC1}) is regular everywhere
(and computable at all zeros of $\sin \pi s$ as well, in preparation to Sect.~\ref{SV2}).
\textsf{Its regularity alone implies that $\ZA (s \sep t)$ 
has the same singular structure as $-\ZS (s \sep t) :$
i.e., only the simple pole $s=1$, of residue~$-\hf$ (property~\tcA)}.

\subsubsection{Special values of $\ZA (s \sep t)$}
\label{SV2}

Once we can compute $[\sin \pi s \mathscr J (s \sep t) ]$ at all integers $s$,
the corresponding values of $\ZA (s \sep t)$ plus $\partial _s \ZA (0 \sep t)$
readily follow from the \tcF-identity (\ref{ZC1}) (as extended to all $s$). 
Then, in the resulting Table~\ref{T3} of special values, the superzeta function $\ZA $
proves fully on par with the much more elementary Hurwitz zeta function of (\ref{ZZD})~!

\textsf{Since $\mathscr J (s \sep t) $ is regular for $s<1$, its transcendental values
get killed in (\ref{ZC1}) by the $\sin \pi s$ factor at all $s=0,-1,-2,\ldots $, 
making those $\ZA (-n \sep t)$ values \textbf{rational} (\tcR)}.

\textit{The remaining, \textbf{transcendental} (\tcT), values $\ZA (n \sep t)$ 
($n=1,2,\ldots$; more detail about $n=1$ in \emph{\cite[\S \,7.4]{VB}})
also directly stem from a generating function
that goes back to a symmetric Hadamard product for} $\Xi (x)$ \cite[\S \,1.10]{E},
\beq
\label{XID}
\Xi (x) = \prod _\rho (1 - x/\rho) \quad \Longrightarrow \quad
\frac{\d}{\d y} \log \Xi (\hf + t + y) \equiv \sum_{n=1}^\infty \ZA (n \sep t) \, (-y)^{n-1} .
\eeq

\begin{table}[h]
\centering
\begin{tabular}  {cccc}
\hline \\[-13pt]
$s$ & $ \ZA (s \sep t) = \sum\limits_\rho (\hf +t -\rho)^{-s} $ & &
\color{dgray} $ \z (s,w) = \sum\limits_{k=0}^\infty (k+w)^{-s} $ \\[7pt]
\hline \\[-11pt]
$-n<0$ & 
$\mathsf{ \frac{\textstyle 2^n}{\textstyle n+1} \, B_{n+1}(\qt + \hf t) + (t + \hf)^n + (t - \hf)^n }$ & &
\color{dgray} $\mathsf{ -\frac{\textstyle 1}{\textstyle n + 1} B_{n+1}(w) }$ \\[5pt]
$0$ & $\mathsf{ \hf(t+\frac{7}{2}) }$ & \tcR & \color{dgray} $\mathsf{ \hf - w }$ \\[3pt]
\hline \\[-13pt]
\hline \\[-12pt]
$\textstyle 0 \atop \mbox{\it \small ($s$-derivative)}$ &
\fbox{$ -\hf (\log 2\pi) \, t + \qt \log 8\pi -\log \Xi(\hf +t) $} & \tcT &
\color{dgray} \fbox{$ -\hf \log 2\pi + \log \G (w) $} \\[6pt]
$\textstyle +1~~ \atop \mbox{\it \small (finite part)}$ &
\fbox{$ \hf \log 2\pi + (\log \Xi)'(\hf +t) $} & & \color{dgray} \fbox{$ - ( \log \G )' (w)$} \\[6pt]
$+1~~$ & $ (\log \Xi)'(\hf +t) $ & & \color{dgray} $\infty$ \\[2pt]
$+n>1$ & $ \displaystyle \frac{(-1)^{n-1}}{(n-1)!} \, (\log \Xi)^{(n)}(\hf +t) $ & &
\color{dgray} $ \displaystyle \frac{(-1)^n}{(n-1)!} \, \bigl( \log \G \bigr) ^{(n)}(w)$\\[4pt]
\hline
\end{tabular}
\caption{\label{T3} \small
Central column: the special values for the superzeta function 
of the 1\textsuperscript{st} kind $\ZA (s \sep t)$,
$ \frac{\mbox{\textsf{Rational}}}{\overline{\mbox{\textit{Transcendental}}}} $ 
at all integer $s=\mp n$.
Right column, for comparison: the same for the Hurwitz zeta function, rewritten from Table~1.
($\mathsf{ B_n(\cdot ): }$ \textsf{Bernoulli polynomials};
$\Xi (\cdot ) :$ \textit{the completed Riemann zeta function (\ref{CZ})}.)
}
\vskip -3mm
\end{table}

The special values $\ZA (n \sep t)$ also display \emph{some imprints} of the fundamental symmetry 
of the Riemann zeros ($\rho \longleftrightarrow (1-\rho), \ t \longleftrightarrow -t$)
(vs \emph{nothing} we know of for non-integer~$s$):

\ni - obvious: \ $\ZA (n \sep {-t}) \equiv (-1)^n \ZA (n \sep {t})$ for $n=1,2, \ldots$,
implying $\ \ZA (n \sep 0) \equiv 0$ for $n$ odd

\ni (vs $\e^{\mi\pi s/2} \!\ZA (s \sep t) - \e^{-\mi\pi s/2} \!\ZA (s \sep {-t})
\equiv 2 \mi (\sin \pi s) \ZC (s \sep \mi t) \ \not\equiv 0$ for non-integer $s$, cf. (\ref{ZCA}));
\smallskip

\ni - less obvious: 
\ $\ZA (-n \sep {-t}) \equiv (-1)^n \ZA (-n \sep {t}) + 
\frac{\textstyle 1}{\textstyle n+1} B_{n+1}(\hf - t) $ 
\quad for ${n=0,1,2, \ldots}$ (as computed from Table~\ref{T3}-\tcR, 
with $B_{n+1} (\cdot) :$ Bernoulli polynomial);
\smallskip

\ni - even less obvious: 
$ 0 \equiv \sum\limits_{k=n}^\infty 
\Bigl( {\textstyle k \!-\! 1 \atop \textstyle n \!-\! 1} \Bigr)
\, t^{k-n} \, \ZA (k \sep t)$ \quad for each odd $n=1,3,\dots$

\ni (skeleton of proof: expand identity $\displaystyle 
(\hf - \rho)^{-n} \equiv (\hf +t - \rho)^{-n} \, [ 1 - t/(\hf +t - \rho) ]^{-n}$ 
in powers of $t$, then sum over $\{\rho \}$, and finally use $\ZA (n \sep 0) \equiv 0$ for $n$ odd).
\smallskip

The next two Tables exemplify Table \ref{T3} for $\ZA _0 (s) \defi \ZA (s \sep 0)$, 
resp. $\ZA _\ast (s) \defi \ZA (s \sep \hf )$.

\begin{table}[h]
\centering
\begin{tabular} {cccl}
\hline \\[-11pt]
$s$ & & $\hfill \ZA _0 (s) \equiv \sum\limits_\rho (\rho-\hf)^{-s} \qquad \quad [t=0] $ & \\[7pt]
\hline \\[-10pt]
& $\!\!\!\!\!\!\!\!\!\!\! $even & $\mathsf{  2^{-n+1}(1-\frac{1}{8} E_n) }$ & \\[-12pt]
$-n \le 0 \ \Biggl\{ $ & & & \\[-15pt]
& $\!\!\!\!\!\!\!\!\!\!\! $odd~ & 
$\mathsf{  -\hf (1 \!-\! 2^{-n}) \frac{\textstyle B_{n+1}}{\textstyle n \!+\! 1} }$ & \\[4pt]
$0$ & & $\mathsf{7/4}$ & \tcR \\[1pt]
\hline \\[-13pt]
\hline \\[-12pt]
$\textstyle 0 \atop \mbox{\it ($s$-derivative)}$ & &
\fbox{$ \ZA '_0 (0) = \log \, \bigl[ 2^{11/4} \pi^{1/2} \G (\qt)^{-1} |\z (\hf)|^{-1} \bigr] $} 
& \tcT \\[6pt]
$\textstyle +1~~ \atop \mbox{\it (finite part)}$ & &
\fbox{$ {\rm FP}_{s=1} \ZA _0 (s) = \hf \log 2 \pi $} & \\[9pt]
& $\!\!\!\!\!\!\!\!\!\!\! $odd~ & $ 0$ & \\[-13pt]
$+n \ge 1 \, \Biggl\{ $ & & & \\[-19pt]
& $\!\!\!\!\!\!\!\!\!\!\! $even & $\!\! \displaystyle 2^{n+1} - \hf \bigl[ (2^n \!-\! 1) \,\z (n) + 2^n \beta (n) \bigr]
- \frac{ (\log |\z |)^{(n)} (\hf) }{(n-1)!}$ & 
$\biggl[ \z (n) \equiv \frac{\textstyle (2\pi)^n |B_n|}{\textstyle 2 \, n!} \biggr]$ \\[8pt]
\hline
\end{tabular}
\caption{\label{T4} \small
$\mathsf{ ( \ {\textstyle E_n \atop \textstyle B_n} \! :
{\textstyle \sf Euler \atop \textstyle Bernoulli \strut} }$ \textsf{numbers}; \quad
$\beta (x) \defi \!\sum\limits_{k=0}^\infty (-1)^k (2k + 1)^{-x} :$ 
\textit{Dirichlet $\beta $-function}.)
}
\end{table}

\begin{table}[h]
\centering
\begin{tabular} {ccc}
\hline \\[-11pt]
$s$ &
$\hfill \ZA _\ast (s) \equiv \sum\limits_\rho \rho^{-s} 
\qquad \qquad \qquad \qquad \ [t=\hf ] $ & \\[7pt]
\hline \\[-10pt]
$-n < 0$ & $\mathsf{  1 - (2^n \!-\! 1) \frac{\textstyle B_{n+1}}{\textstyle n \!+\! 1} }$ & \\[3pt]
$0$ & $\mathsf{2}$ & \tcR \\[1pt]
\hline \\[-13pt]
\hline \\[-11pt]
$\textstyle 0 \atop \mbox{\it ($s$-derivative)}$ & \fbox{$\ZA' _\ast (0) = \hf \log 2 $} & \tcT \\[6pt]
$\textstyle +1~~ \atop \mbox{\it (finite part)}$ &
\fbox{$ {\rm FP}_{s=1} \ZA _\ast (s) = 1 - \hf \log 2 + \hf \g $} & \\[6pt]
$+1~~$ & $  1 -\hf \log 4\pi + \hf \g $ & \\
$+n > 1$ &
$ \displaystyle \!\!\! 1 - (1 \!-\! 2^{-n})\,\z (n) + \frac{g^{\rm c}_n}{(n - 1)!} \ \equiv \
1 - (-1)^n 2^{-n} \z (n) - \frac{ (\log |\z |)^{(n)}(0) }{(n-1)!}$ & \\[8pt]
\hline
\end{tabular}
\caption{\label{T5} \small
($\g :$ \textit{Euler's constant; 
\ $g^{\rm c}_n :$ cumulants of the Stieltjes constants} 
\cite[\S \,5]{Is}\cite{BL}\cite[\S \,3.3]{Vz}\cite{C2}, \textit{of generating function
$\ \log \, [y \, \z (1+y)] \equiv {- \! \sum\limits_{n=1}^\infty 
\! \frac{\textstyle (-1)^n}{\textstyle n!} g_n^{\rm c} \, y^n } \,$; 
\ $\z (n) \equiv \frac{\textstyle (2\pi)^n |B_n|}{\textstyle 2 \, n!}$ for $n$ even}.)
}
\end{table}

\subsection{The superzeta function of the 2\textsuperscript{nd} kind $\ZB (\s \sep t)$}

Defined by (\ref{ZBD}) for $\Re \s > \hf$, $\ZB (\s \sep t)$ appears (for $t \ne 0$) 
functionally independent from the function of the 1\textsuperscript{st} kind $\ZA (s \sep t)$;
contrary to the latter, it manifestly embodies the full symmetry
of the Riemann zeros, through the identity $\ZB (\s \sep t) \equiv\ZB (\s \sep {-t})$.

On the other hand, $\ZB (\s \sep t)$ for $t \ne 0$ displays no tractable functional relation~\tcF.
Still, by expansion in powers of $t$ around the known confluent case (\ref{Z0}) 
\cite[\S \,6.1]{Vz}\cite[\S \,8.3]{VB}, $\ZB (\s \sep t)$ \textsf{is proved meromorphic in~$\s $, 
now with \emph{double} poles $\s = \hf, -\hf, -\frac{3}{2}, \ldots$, 
and all principal parts plus rational values 
$\ZB (-m \sep t) \ (m \in {\mathbb N})$ explicitly computable (\tcA+\tcR).
E.g., the leading pole $\s = \hf $ has the principal part ($t$-independent) }
\beq
\label{POL}
\ZB (\hf + \eps \sep t) =
\frac{1}{8 \pi} \ \eps^{-2} -\frac{\log 2 \pi}{4 \pi} \ \eps^{-1} \quad 
\mbox{\textsf{+ [regular term]}}, \quad \eps \to 0 .  \qquad \qquad \mbox{(\tcA)}
\eeq

\subsubsection{Special values of $\ZB (\s \sep t)$ (continued)}

\itshape
The remaining, \textbf{transcendental} (\tcT), values $\ZB (m \sep t)\ (m=1,2,\ldots)$ 
stem from a variant of~(\ref{XID}) \emph{\cite[\S \,1.10]{E}},
\beq
\label{XIE}
\Xi (x) = \!\!\prod_{\Im \! \rho >0} \Bigl[ 1 - \frac{x(1-x)}{\rho (1-\rho)} \Bigr] \quad \Rightarrow \quad
\frac{\d}{\d w} \log \Xi (\hf + (t^2+w)^{1/2}) \equiv \sum_{m=1}^\infty \ZB (m \sep t) \, (-w)^{m-1} ;
\eeq
only the derivative $\partial_\s \ZB (\s \sep t) |_{\s =0}$ is harder to obtain
(see \emph{\cite[\S \,4]{Vz}\cite[\S \,8.4]{VB}}), and that completes Table~\ref{T6}.

The change of variables $t+y = (t^2+w)^{1/2}$ in the generating function (\ref{XID}) 
(at fixed $t$, say $t \ge 0$) also has to yield (\ref{XIE}), 
hence the $\ZB (m \sep t) \ (m=1,2,\ldots)$ must relate to the $\ZA (n \sep t) \ (1 \le n \le m)$.
To get a fully explicit connection, our shortcut (new) is to use the residue calculus twice, 
inside small positive contours encircling $y=0$, resp. $w=0 :$
\beq
\label{ZAR}
\frac{(-1)^{n-1}}{2\pi\mi} \oint \frac{\d \log \Xi (\hf + t + y)}{y^n} = \biggl\{ \matrix{
\ZA (n \sep t) & \mbox{for } n=1,2,\ldots , \hfill \cr 0 & \mbox{for $n \in \mathbb Z$ otherwise} ; }
\eeq
\bea
\ZB (m \sep t) \si=\se \frac{(-1)^{m-1}}{2\pi\mi} \oint \frac{\d \log \Xi (\hf + (t^2+w)^{1/2})}{w^m}
= \frac{(-1)^{m-1}}{2\pi\mi} \oint \frac{\d \log \Xi (\hf + t + y)}{(2ty + y^2)^m} \nonumber \\
\si=\se \frac{(-1)^{m-1}}{2\pi\mi} \! \oint \!
\sum_{\tilde n=-m}^{\infty} \! (-1)^{m+\tilde n} {2m \!+\! \tilde n \!-\! 1 \choose m-1} 
(2t)^{-2m-\tilde n} y^{\tilde n} \ \d \log \Xi (\hf \!+\! t \!+\! y) \quad (t \ne 0) \nonumber \\
\label{ZBA}
\si\equiv\se \sum_{n=1}^m { 2m \!-\! n \!-\! 1 \choose m-1} (2t)^{-2m+n} \! \ZA (n \sep t) 
\quad \mbox{for } m=1,2,\ldots , \mbox{ and } \ t \ne 0
\eea
by (\ref{ZAR}) (vs (\ref{Z1E}) below for $t=0$); e.g., the first of the identities (\ref{ZBA}) is
$\ZB (1 \sep t ) \equiv \frac{\textstyle \ZA (1 \sep t)}{\textstyle 2t \strut}$.

\normalfont

\begin{table}[h]
\centering
\begin{tabular}  {ccc}
\hline \\[-13pt]
$\sigma$ &
$ \ZB (\s \sep t) = \sum\limits_{k=1}^\infty ({\tau_k}^2+ t^2)^{-\s } $ & \\[8pt]
\hline \\[-12pt]
$-m<0$ & $\mathsf{  (t^2 - \qt)^m 
- 2^{-2m-3} \sum\limits_{j=0}^m  {\textstyle m \vphantom{j} \choose \textstyle j} 
(-1)^j \, E_{2j} \, (2t)^{2(m-j)} }$ & \\
$0$ & $\mathsf{7/8}$ & \tcR \\
\hline \\[-13pt]
\hline \\[-12pt]
$\textstyle 0 \atop \mbox{\it ($\s $-derivative)}$ &
\fbox{$ \partial _\s \ZB (0 \sep t) = \qt \log 8\pi -\log \Xi(\hf +t) $}  & \tcT \\[6pt]
$+m \ge 1$ & 
$\displaystyle \frac{(-1)^{m-1}}{(m-1)!} \, \frac{\d^m}{\d(t^2)^m} (\log \Xi )(\hf +t) $ & \\[8pt]
\hline
\end{tabular}
\caption{\label{T6} \small
The special values for the superzeta function 
of the 2\textsuperscript{nd} kind $\ZB (\s \sep t)$,
$ \frac{\mbox{\textsf{Rational}}}{\overline{\mbox{\textit{Transcendental}}}} $ 
at all integer $\s=\mp m$.
($\mathsf{ E_n : }$ \textsf{Euler numbers};
$\Xi (\cdot ) :$ \textit{the completed Riemann zeta function (\ref{CZ})}.)
}
\end{table}

All in all and notwithstanding its lack of a tractable functional relation \tcF, 
the superzeta function of the 2\textsuperscript{nd} kind matches 
its partner of the 1\textsuperscript{st} kind 
(or the Hurwitz zeta function) for its stock of explicit fixed-$\s $ properties.
\medskip

Next, we exemplify those formulae at the special parameter values $t=0$ and $\hf$.

For $\ZB _0 (\s) = \ZB (\s \sep t=0)$ (the confluent case), Table \ref{T4} can be reused 
since (\ref{Z0}) entails
\beq
\label{Z1E}
\ZB _0 (m) \equiv \hf (-1)^m \ZA _0 (2m) \quad (m \in {\mathbb Z}), \qquad \qquad 
\frac{\d}{\d \s } \ZB _0 \, (0) = \frac{\d}{\d s} \ZA _0 \, (0) .
\eeq

Whereas for $\ZB _\ast (\s) = \ZB (\s \sep t=\hf )$, Table~\ref{T6} plus (\ref{ZBA}) specialize to:

\begin{table}[h]
\centering
\begin{tabular}  {cccc}
\hline \\[-13pt]
$\sigma$ & & \qquad \qquad \quad 
$\ZB _\ast(\s ) = \sum\limits_{k=1}^\infty ({\tau_k}^2+ \qt)^{-\s } \qquad [t=\hf]$ & \\[8pt]
\hline \\[-12pt]
$-m < 0$ & &
$\mathsf{  -2^{-2m-3} \sum\limits_{j=0}^m  
{\textstyle m \vphantom{j} \choose \textstyle j} (-1)^j \, E_{2j} }$ & \\[6pt]
$0$ & & $\mathsf{7/8}$ & \tcR \\
\hline \\[-13pt]
\hline \\[-11pt]
$\textstyle 0 \atop \mbox{\it ($\s $-derivative)}$ & & 
\fbox{$ \ZB '_\ast (0) = \qt \log 8\pi $} & \tcT \\[6pt]
$+m \ge 1$ & & 
$ \sum\limits_{n=1}^m \Bigl( {\textstyle 2m \!-\! n \!-\! 1 \vphantom{j} \atop \textstyle m-1} \Bigr)
\ZA _\ast (n)$ & \\[8pt]
\hline
\end{tabular}
\caption{\label{T7} \small
As Table \ref{T6}, but at $t=\hf$.
($\mathsf{ E_n : }$ \textsf{Euler numbers}; \textit{values $\ZA _\ast(n) :$ see Table~\ref{T5}.})
}
\end{table}

\subsection{The superzeta function of the 3\textsuperscript{rd} kind $\ZC (s \sep \tau)$}

Defined by (\ref{ZCD}) for $\Re s > 1$, $\ZC (s \sep \tau)$ is less regular 
than the other two kinds \cite[\S \,3]{HKW}\cite[\S \,6.2]{Vz},
nevertheless it reduces to functions of the 1\textsuperscript{st} kind: 
indeed, the straightforward identity
\beq
\ZA (s \sep t) \equiv \e^{\mi\pi s/2} \ZC ( s \sep \mi t) + \e^{-\mi\pi s/2} \ZC (s \sep {-\mi t}) 
\eeq
can be inverted, yielding \cite[\S \,8]{Me}\cite[Chap.~9]{VB}
\beq
\label{ZCA}
\ZC (s \sep \tau) \equiv \frac{1}{2 \mi \sin \pi s} \,
\bigl[ \e^{\mi\pi s/2} \ZA (s \sep {-\mi \tau}) - \e^{-\mi\pi s/2} \ZA (s \sep {+\mi \tau}) \bigr] .
\eeq
The singular structure of $\ZC $ is thus computable just as before;
on the other hand, only one special value remains explicit: 
the finite part of $\ZC (s \sep \tau)$ at $s=0$, which reads
$\displaystyle \frac{7}{8} + \frac{\log 2 \pi}{2 \pi} \, \tau $.

\subsection{Extensions}
\label{SE}

Those currently include, but should not be limited to, (zeta functions over the) zeros of:
\bigskip

\ni \textbullet \ Dedekind zeta functions (for algebraic number fields) 
\cite{Ku}\cite{HKW}\cite{HIKW}\cite{Ih}\cite{VO}
\medskip

\ni \textbullet \ Dirichlet L-functions (for real primitive Dirichlet characters) \cite{VO}
\medskip

\ni \textbullet \ Selberg zeta functions (for cocompact subgroups of SL(2,$\mathbb R$))
- with more work (\cite[App.~B]{VB} and refs. therein).
\medskip

In every such extension, $\ZA _\ast (1) = \sum_\rho 1/\rho $ 
is closely related to the associated \emph{generalized Euler constant}
(= \emph{Euler--Kronecker invariant}), a major invariant of the zeta function \cite{HIKW}\cite{Ih}.

\vskip 1cm

\section{Complex-asymptotic view of Riemann Hypothesis}

\subsection{The Keiper--Li coefficients $\l _n$}

The sequence of real numbers
\beq
\l _n = \sum_\rho \bigl[ 1-(1-1/\rho)^n \bigr] \qquad \qquad \qquad (n=1,2,\ldots) ,
\eeq
(in Li's notation \cite{Li1}, amounting to $n$ times Keiper's $\l _n$ \cite{Ke}), 
of generating function
\beq
\label{LId}
\log \Xi \Bigl( x= \frac{1}{1-z}\Bigr) \equiv \sum_{n=1}^\infty \frac{\l _n}{n} \, z^n ,
\eeq
has served to recast the Riemann Hypothesis (RH) (cf. (\ref{RRH})), as

\centerline{RH $ \iff \ \l _n >0 \ $ for all $n$ \qquad (\emph{Li's criterion}) \cite{Li1}\cite{BL}.} 

Here we review the (large-$n$) \emph{asymptotics} of the $\l _n$ instead \cite{Ke}\cite{M1}\cite{VL}\cite{Lg},
and our resulting \emph{asymptotic criterion for RH} that involves just the \emph{tail} of the sequence $\{ \l _n \}$
\cite{VL}\cite[Chap.~11]{VB}.

\subsection{Exponential-asymptotic large-$n$ analysis of $\l _n$}
\label{EAA}

\begin{figure}[h]
\centering
\includegraphics[scale=.5]{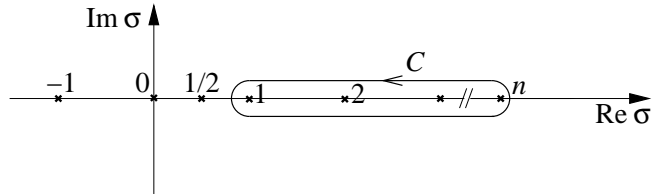}
\caption{\label{F3} \small
Contour of integration for (\ref{LIN}). 
( \protect\raisebox{1pt}{$\scriptscriptstyle \pmb{\times}$} : poles of the integrand.)
}
\end{figure}

The Hadamard product formula in (\ref{XIE}) implies
\vskip -2mm

\beq
\label{XIZ}
\log \Xi \Bigl( {1 \over 1-z} \Bigr) =
\sum_{k=1}^{\infty} \log \biggl[ 1 + \frac{z }{({\tau _k}^2 \!+\! \qt ) (1-z)^2 } \biggr] =
- \sum_{m=1}^{\infty} \frac{\ZB _\ast(m)}{m} \frac{(-z)^m}{(1-z)^{2m}} 
\eeq
(the $\ZB _\ast (m)$ are in Table~\ref{T7}); then, expanding $(1-z)^{-2m}$ by the generalized binomial formula,
reordering in powers of $z$ and substituting into (\ref{LId}), we get \cite{VL}\cite[\S \,8.6.2]{VB}
\vskip -2mm

\beq
\label{LZ2}
\l _n = - n \sum_{m=1}^n {(-1)^m \over m}{m+n-1 \choose 2m-1} \, \ZB _\ast(m) , \qquad n=1,2,\ldots .
\eeq
This yields a representation by an integral over the superzeta function 
$\ZB _\ast (\s ) \equiv \ZB (\s \sep \hf )$ (of the 2\textsuperscript{nd} kind),
\vskip -3mm

\beq
\label{LIN}
\l _n = \frac{(-1)^n n \,\mi}{\pi} 
\oint_C \frac{ \G (\s +n) \G (\s -n)}{\G (2 \s +1) } \, \ZB _\ast(\s ) \,\d \s ,
\eeq
as proved by reduction to (\ref{LZ2}) using the residue calculus, see Fig. \ref{F3}.
Our point is now that for $n \to +\infty$ (and after using the Stirling formula in the numerator),
the integral form (\ref{LIN}) is asymptotically computable 
thanks to the \emph{method of steepest descents} \cite{Er}. 

With the integrand in (\ref{LIN}) behaving like $n^{2 \s -1}$ for $n \to +\infty$,
that means to push the integration path toward decreasing $\Re \s $, up to \emph{saddle-points} 
of the integrand where dominant contributions to the integral will localize.
Here, saddle-points $\s _\ast$ of two species compete in the half-plane $\{ \Re \sigma >\hf \}$ 
for dominance (= largest $\Re \s _\ast $), see Fig.~\ref{F4}~\cite{VL}:

\begin{figure}[h]
\centering
\includegraphics[scale=.5]{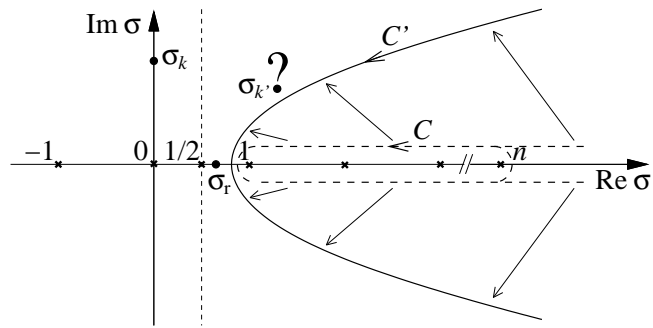}
\caption{\label{F4} \small
Integration-path deformation for a large-$n$ evaluation of the integral~(\ref{LIN}).
\protect\raisebox{1pt}{$\ \scriptstyle \bullet $} : typical saddle-points (none on scale); 
$\s _k (n)$ comes from a Riemann zero $\hf +\mi \tau _k$ on the critical line 
(and is irrelevant, beyond reach),
$\s _{k'} (n)$ comes from a \emph{putative} (``?") zero $\hf +\mi \tau _{k'}$ 
off the critical line (and then, gives a dominant contribution),
finally $\sigma_{\rm r}(n)$ is the real saddle-point tending to the pole $\s = 1/2 $ from above
(and gives the dominant contribution if and only if RH holds).
}
\end{figure}

\ni \textbullet \ using $\ZB _\ast(\sigma) = \sum_k ({\tau_k}^2+ \qt)^{-\s }$ 
to integrate (\ref{LIN}) term by term:
then the $k$-th integrand has a saddle-point $\sigma_k(n) \sim \hf \, n \, \mi / \tau_k $,
which is relevant iff $\arg \tau_k >0$, as this is the necessary and sufficient condition
for that saddle-point to ultimately reside in the half-plane $\{ \Re \sigma >\hf \}$  
(that is, for $n \gtrsim |\Im 1/\tau_k| ^{-1}$);
\medskip

\ni \textbullet \ the pair of poles $\s =+\hf$  of $\ZB _\ast $ as in (\ref{POL}) 
and $\s =1$ of $\G (\s -n)$
together force a \emph{real} local minimum in the modulus of the integrand, unconditionally:
$ \s _{\rm r}(n) \sim \hf + \frac{1}{\log n} $, 
this is another saddle-point, and the dominant one iff all $\arg \tau_k \equiv 0$ (i.e., RH true).
\smallskip

Thus, the competition between both types of saddle-point contributions to (\ref{LIN}) 
results in this \emph{asymptotic criterion for RH} \cite{VL}: as $n \to +\infty$,
\bea
\mbox{- \bf if RH is false,} \si \l _n \se 
\mbox{has an \emph{exponentially growing} oscillatory behavior:} \nonumber  \\
\label{RS}
\l  _n \si\sim\se - \!\! \sum_{\{ \arg \tau_k >0 \}}
\Bigl( \frac{\tau_k+\mi/2}{\tau_k-\mi/2} \Bigr) ^n
\! + \mbox{c. c.} \qquad \pmod {{\rm o}(\e^{\eps n}) \ \forall \eps >0} \\[3pt]
\mbox{- \bf if RH is true,} \si \l _n \se \mbox{has a \emph{tempered} growth to $+\infty :$} \nonumber  \\
\label{RER}
\l  _n \si\sim\se \hf n \, (\log n - 1+\g   - \log 2\pi) \qquad \pmod {{\rm o}(n)} .
\eea

However, if a numerical crossover from the latter (RH) to the former (non-RH) behavior
is sought as a signal of a violation of RH,
that has to await a huge value $n \approx 10^{18}$ at least,
corresponding to the present height up to which $\Re \rho = \hf $ has been verified \cite{Lg}.

\subsection{Extensions}
\label{LIX}

For zeros of a more general zeta function as in Sect. \ref{SE}, and such that
\beq
\ZB _\ast (\hf+\eps) = R_{-2} \,\eps^{-2} + R_{-1} \,\eps^{-1} + O(1)_{\eps \to 0} \qquad
(\mbox{with }R_{-2}>0),
\eeq
(generalizing (\ref{POL}) from the case of the Riemann zeros),
then the \emph{Generalized Riemann Hypothesis} is equivalent to
\bea
\l _n &\sim& (-1)^n \, 2n \, 
\Res_{\s =1/2} \Biggl[ \frac{\G (\s +n)\G (\s -n)}{\G (2\s +1)} \, \ZB _\ast (\s ) \Biggr] \\ 
\nonumber \\
&\sim& 2 \pi n \, [ 2 R_{-2} (\log n - 1 + \gamma ) + R_{-1} ] \qquad \qquad \qquad \qquad \pmod{o(n)} 
\eea
for its corresponding Keiper--Li sequence $\{\l _n\}$ as $n  \to \infty$ \cite{Lg};
the constants $R_{-2},\ R_{-1}$ are also those governing the asymptotic \emph{counting function} 
of the zeros' ordinates, as
\beq
N(T) \sim 2T \, [ 2 R_{-2} (\log T -1) + R_{-1} ] \qquad (T \to +\infty).
\eeq

\subsection{(New) Shifted asymptotic and Li criteria for RH}

\ni (returning to the Riemann-zeta case alone, for the sake of definiteness).

The Keiper--Li (KL) generating series (\ref{LId}) combines two features:
a high sensitivity to the Riemann Hypothesis through its radius of convergence
($|z|=1$ if and only if RH holds vs \emph{strictly less} otherwise),
obtained by the use of a conformal mapping $z \mapsto x$ that sends 
the disk $\{ |z|<1 \}$ to the half-plane $\{ \Re x > \hf \}$;
\emph{and} an expansion basepoint (the image of $z=0$) set to $x=1$
(hence the parameter $t \equiv x-\hf$ is $\hf$ in the superzeta function entering \S\,\ref{EAA}: 
$\ZB _\ast (\s ) \defi \ZB (\s \sep \hf )$). The two features always appeared coupled,
as the conformal mapping was invariably prescribed to be $x=1/(1-z)$, as in (\ref{LId}).

Yet, we claim that imposing only $x=1$ as basepoint is an \emph{unnecessary limitation}.
To wit, we present an alternative KL-like framework based at the other remarkable point 
$x=\hf $ (the center of symmetry), i.e., $t=0$. In the last minute, we saw
that the idea of a general basepoint had been developed elsewhere, only recently 
and in a different manner that just excludes our preferred symmetrical point $x=\hf $ \cite{Se}.

For expanding $\Xi (x)$ about $x=\hf $, due to (\ref{RFE}) the good variable is $w=t^2$,
in which the above half-plane $\{ \Re x > \hf \}$ maps to the cut $w$-plane 
${\mathbb C} \setminus {\mathbb R}^-$; 
this cut can then be truncated to any half-line $(-\infty,-w']$ 
with $0 < w' < \inf_\rho \{ |\Im \rho |^2 \}$ 
(that only deletes a zero-free segment), and here we choose $w'=\qt $
just to make the formulae simpler and closest to the KL case. 
There now remains to map to the cut $w$-plane ${\mathbb C} \setminus (-\infty,-\qt ] $ 
from a unit disk $\{ |y|<1 \}$, which is achieved by $w=y/(1-y)^2$. 
Under the resulting $y \mapsto x$ correspondence,
we will need the preimages of a Riemann zero $\rho =\hf + \mi \tau _k$, as
\vskip -2mm

\beq
\label{YPM}
y_{k,\pm} = 1 - \frac{1}{2 {\tau_k}^2} \pm \frac{\mi}{\tau _k} 
\Bigl( 1 - \frac{1}{4 {\tau_k}^2} \Bigr) ^{1/2} \quad \equiv \e^{\pm \mi \theta_k} 
\ \mbox{ with } \theta _k \defi 2 \arcsin \frac{1}{2 \tau _k} .
\eeq

Now we start from the generating function (\ref{XIE}) at $t=0$ and integrate in $w$, to get
$ \log \Xi _0 (\hf +w^{1/2}) \equiv 
- \sum\limits_ {m=1}^\infty \frac{\textstyle \ZB _0(m)}{\textstyle m^{\vphantom{o}}} (-w)^m $
where $\Xi _0(x) \defi \Xi (x)/\Xi (\hf )$, amounting to
\vskip -.5mm

\beq
\label{X0Z}
\sum_{n=1}^\infty \frac{\l _n^0}{n} \, y^n \defi
\ \ \log \Xi _0  \Bigl( \frac{1}{2} + \frac{y^{1/2}}{1-y} \Bigr) \equiv
- \sum_ {m=1}^\infty \frac{\ZB _0(m)}{m} \frac{(-y)^m}{(1-y)^{2m}} \, ;
\eeq
thus, we created a KL-like sequence of ``central" coefficients $\l _n^0$. 
Comparing to (\ref{LId})--(\ref{XIZ}),
we see that the relationships (\ref{LZ2}) persist with the $\l _n$ substituted by the $\l _n^0$ 
and the superzeta function $\ZB _\ast$ by $\ZB _0$.
(Equivalently, by (\ref{Z1E}) and Table~\ref{T3} at $t=0$,
\beq
\label{LZ1}
\qquad \l _n^0 = - \frac{n}{2} \sum_{m=1}^n {m\!+\! n \!-\! 1 \choose 2m-1} \frac{1}{m} \ZA _0(2m) 
= n  \sum_{m=1}^n {m\!+\! n \!-\! 1 \choose 2m-1}  \frac{1}{(2m)!} \, (\log \Xi )^{(2m)} (\hf) ,
\eeq
and the last line in Table~\ref{T4} displays $\ZA _0(2m)$ in most reduced form.) E.g.,
\bea
\l _1^0 \si=\se \ZB _0(1) = \hf (\log \Xi )''(\hf) \approx 0.0231050, \\
\l _2^0 \si=\se 4 \ZB _0(1) - \ZB _0(2) = 2 (\log \Xi )''(\hf) 
+ {\textstyle \frac{1}{12}} (\log \Xi )^{(4)}(\hf) \approx 0.0923828 , \ \ldots 
\eea

Being based upon (\ref{LZ2}), our asymptotic analysis of $\{ \l _n \}$ 
sketched in \S\,\ref{EAA} \cite{VL}\cite[Chap.~11]{VB}
then carries over likewise to $\{ \l _n^0 \}$, 
and it now yields this \emph{shifted asymptotic criterion} for RH:
\vskip -5mm

\bea
\label{RS0}
\si \mbox{\bf -} \se \mbox{\textbf{if RH is false,}} \quad 
\l _n^0 \sim - \!\! \sum_{\{ \arg \tau_k >0 \}} y_{k,-}^{\ -n} + \, \mbox{c. c.} , \\
\noalign{\ni an exponentially growing oscillatory behavior (like (\ref{RS}))
since $|y_{k,-}^{\, -1}|>1$ if $\arg \tau_k >0$;} 
\nonumber\\[-8pt]
\si \mbox{\bf -} \se \mbox{\textbf{if RH is true,}} \quad \ \l _n^0 \sim \hf n \, (\log n - 1+\g   - \log 2\pi) ,
\eea
the latter being unchanged from (\ref{RER}), as according to \S\,\ref{LIX} it arises from 
the principal part (\ref{POL}) of $\ZB (\s \sep t)$ at $\s =\hf $, and (\ref{POL}) is $t$-independent.
\smallskip

Finally, \emph{Li's criterion} \cite{Li1} shifts as well: by the Hadamard product formula in (\ref{XIE}),
\bea
\log \Xi _0 \Bigl( \frac{1}{2} + \frac{y^{1/2}}{1-y} \Bigr) \si=\se 
\sum_k \log \biggl[ 1 + \frac{y}{{\tau_k}^2 (1-y)^2} \biggr] \\
\si=\se \sum_k \bigl( -\log (1-y)^2 + \log \, \bigl[ (1-y/y_{k,-})(1-y/y_{k,+}) \bigr] \bigr) \\[-2pt]
\si=\se \sum_{n=1}^\infty \Biggl[ \sum_k  (2 - y_{k,-}^{\ -n} - y_{k,+}^{\ -n}) \biggr] \frac{y^n}{n} \\[4pt]
\label{L0C}
\Longrightarrow \qquad \l _n^0 \si\equiv\se \sum_k  (2 - y_{k,-}^{\ -n} - y_{k,+}^{\ -n}) \quad
\mbox{by identification with (\ref{X0Z}).}
\eea
Now if RH is true, then for all $k$, $\tau _k$ is purely real and so is $\theta _k $ in (\ref{YPM}),
in which case $y_{k,\pm} = \e^{\pm \mi \theta_k}$ makes (\ref{L0C}) manifestly positive for all $n$. 
Whereas if RH is false, then as $n \to +\infty$, the growing oscillatory behavior (\ref{RS0}) 
will force some $\l _n^0$ into the negative range. 
So, the Li criterion transfers to this ``central" sequence:
\beq
\mbox{RH} \qquad \iff \qquad \l _n^0 >0 \quad  \mbox{for all } n .
\eeq

{\bf Remark:}
Other expressions for the $\l _n^0$, derivable from (\ref{X0Z}) 
by using the residue calculus both ways (similarly to the proof of (\ref{ZBA})), are
\bea
\qquad \l _n^0 \si\equiv\se \frac{n}{2 \pi\mi} 
\oint \log \Xi _0 \Bigl( \frac{1}{2} + \frac{y^{1/2}}{1-y} \Bigr) \frac{\d y}{y^{n+1}} 
\qquad \qquad \qquad \qquad \qquad (n=1,2,\ldots)\\
\si\equiv\se \frac{n}{4 \pi\mi} \oint \log \Xi _0 (\hf +t) \, (t^2+\qt )^{-1/2} 
\Bigl[ \frac{t^2}{y(t)} \Bigr] ^n \frac{\d t}{t^{2n+1}} \nonumber\\
\noalign{\qquad (upon the change of variable $\displaystyle t \equiv \frac{y^{1/2}}{1-y}
\mbox{ implying } \d y (t) = (t^2+\qt )^{-1/2} y \, \frac{\d t}{t} $)}
\si\equiv\se \frac{1}{4 \, [2(n-1)]!} \Bigl( \frac{\d}{\d t} \Bigr)^{2n} \biggl[ 
(t^2+\qt )^{-1/2} \bigl( 1 - y(t) \bigr) ^{-2n} \log \Xi _0 (\hf +t) \biggr]_{t=0} \, .
\eea

\section*{Acknowledgements}

We are most grateful to the conference organizers and to
the Research Institute for Mathematical Sciences of Kyoto University
for their kind invitation and generous support.

\end{document}